\documentclass[12pt]{article}
\usepackage{amsfonts}

\input{tcilatex}
\begin{document}

\begin{center}
{\Large On McMillan's theorem }

{\Large about uniquely decipherable codes}

\bigskip

Stephan Foldes

Tampere University of Technology

PL 553, 33101 Tampere,\ Finland

sf@tut.fi

19 June 2008

\bigskip

\bigskip \textbf{Abstract}
\end{center}

\textit{Karush's proof of McMillan's theorem is recast as an argument
involving polynomials with non-commuting indeterminates certain evaluations
of which yield the Kraft sums of codes, proving a strengthened version of
McMillan's theorem.}

\bigskip

\bigskip

Let $len:A^{\ast }\longrightarrow 
\mathbb{N}
$ be the length function on the free monoid of all strings over a given
non-empty finite set $A.$ Let $con:A^{\ast \ast }\longrightarrow A^{\ast }$
be the concatenation map which to every string of strings associates their
concatenation. A \textit{uniquely decipherable code} is a finite set $%
C\subseteq A^{\ast }$ such that $con$ is injective on the submonoid $C^{\ast
}$ of $A^{\ast \ast }$. This submonoid is then isomorphic to the submonoid $%
\overline{C}=con[C^{\ast }]$ of $A^{\ast }$ freely generated by $C.$
Denoting by $r$ the number of elements of the \textit{alphabet} $A$, the 
\textit{Kraft sum} $K(C)$ of any finite $C\subseteq A^{\ast }$ is defined as 
$\sum_{x\in C}$ $r^{-len(x)}.$ In [M] McMillan showed that if $C$ is a
uniquely decipherable code, then its Kraft sum is at most $1.$ The proof
usually given is that of Karush [K]. This proof can be recast as an argument
involving evaluations of polynomials with non-commuting indeterminates
corresponding to the various (infinitely many) strings in $A^{\ast }$, as
follows.

\bigskip

Let $%
\mathbb{R}
\left\langle A^{\ast }\right\rangle $ be the free associative $%
\mathbb{R}
$-algebra generated by the elements of $A^{\ast }$ considered as
indeterminates, i.e.\ $%
\mathbb{R}
\left\langle A^{\ast }\right\rangle $\ is the non-commutative ring of formal
polynomials with real coefficients in the non-commuting indeterminates $x\in
A^{\ast }.$ For $w=(x_{1},...,x_{n})$ in $A^{\ast \ast }$ let $P(w)$ denote
the monomial $x_{1}\cdot ...\cdot x_{n}$ in $%
\mathbb{R}
\left\langle A^{\ast }\right\rangle $.

Let $C$ and $D$ be finite uniquely decipherable codes over a non-empty
finite alphabet $A$ with $r$ elements, and suppose that $C\subseteq 
\overline{D}.$ The Kraft sum $K(C)$ of $C$ is then the evaluation of the
polynomial $\sum_{x\in C}x$ at $x:=r^{-len(x)}$ for $x\in A^{\ast }$.

Fix a positive integer $k$. For any positive integer $l$, partition the set $%
D^{l}$ into two disjoint sets:%
\[
W_{l1}=\left\{ w\in D^{l}:con(w)\in con[C^{k}]\right\} 
\]%
\[
W_{l2}=D^{l}\setminus W_{l1}
\]%
For every $l$ the polynomial 
\[
\left( \tsum\limits_{x\in D}x\right) ^{l}=\tsum\limits_{w\in D^{l}}P(w)\text{
\ \ \ \ \ \ \ \ \ \ \ \ \ \ \ \ \ \ \ \ \ \ \ \ \ \ \ \ \ \ \ \ \ \ \ \ \ \
\ \ \ \ \ \ \ \ (1)}
\]%
is equal to the sum%
\[
\tsum\limits_{w\in W_{l1}}P(w)+\tsum\limits_{w\in W_{l2}}P(w)\text{ \ \ \ \
\ \ \ \ \ \ \ \ \ \ \ \ \ \ \ \ \ \ \ \ \ \ \ \ \ \ \ \ \ \ \ \ \ \ \ \ \ \
\ (2)}
\]%
Let $m$ be the largest integer $n$ with $C\cap con[D^{n}]\neq \emptyset .$
Then the polynomial 
\[
\tsum\limits_{l=k}^{mk}\left( \tsum\limits_{x\in D}x\right)
^{l}=\tsum\limits_{l=k}^{mk}\tsum\limits_{w\in D^{l}}P(w)\text{ \ \ \ \ \ \
\ \ \ \ \ \ \ \ \ \ \ \ \ \ \ \ \ \ \ \ \ \ \ \ \ \ \ \ \ \ (3)}
\]%
is the sum of%
\[
\tsum\limits_{l=k}^{mk}\tsum\limits_{w\in W_{l1}}P(w)\text{ \ \ \ \ \ \ \ \
\ \ \ \ \ \ \ \ \ \ \ \ \ \ \ \ \ \ \ \ \ \ \ \ \ \ \ \ \ \ \ \ \ \ \ \ \ \
\ \ \ \ \ \ \ \ \ \ \ \ (4)}
\]%
and%
\[
\tsum\limits_{l=k}^{mk}\tsum\limits_{w\in W_{l2}}P(w)\text{ \ \ \ \ \ \ \ \
\ \ \ \ \ \ \ \ \ \ \ \ \ \ \ \ \ \ \ \ \ \ \ \ \ \ \ \ \ \ \ \ \ \ \ \ \ \
\ \ \ \ \ \ \ \ \ \ \ \ (5)}
\]%
Let $I(C,D)$ be the ideal of $%
\mathbb{R}
\left\langle A^{\ast }\right\rangle $\ generated by the polynomials $x-P(w)$
for $x\in C,w\in D^{\ast },x=con(w).$ Modulo this ideal, (4) is congruent to 
\[
\left( \tsum\limits_{x\in C}x\right) ^{k}\text{ \ \ \ \ \ \ \ \ \ \ \ \ \ \
\ \ \ \ \ \ \ \ \ \ \ \ \ \ \ \ \ \ \ \ \ \ \ \ \ \ \ \ \ \ \ \ \ \ \ \ \ \
\ \ \ \ \ \ \ \ \ \ \ \ (6)}
\]

\bigskip

The homomorphism $%
\mathbb{R}
\left\langle A^{\ast }\right\rangle \longrightarrow 
\mathbb{R}
$ evaluating each polynomial at $x:=r^{-len(x)}$ is null on the ideal $%
I(C,D) $ and therefore the evaluation of (3) equals the sum of the
evaluations of (5) and (6). The evaluation of (5) being non-negative, the
evaluation of (6) is at most the evaluation of (3). For the Kraft sums $K(C)$
and $K(D)$ this means that%
\[
K(C)^{k}\leq \tsum\limits_{l=k}^{mk}K(D)^{l}\text{ \ \ \ \ \ \ \ \ \ \ \ \ \
\ \ \ \ \ \ \ \ \ \ \ \ \ \ \ \ \ \ \ \ \ \ \ \ \ \ \ \ \ \ \ \ \ \ \ \ \ \
(7)} 
\]%
Applying this to $D=A^{1}$, as $C\subseteq \overline{A^{1}}$ and obviously $%
K(A^{1})=1$, we obtain%
\[
K(C)^{k}\leq mk-k+1\leq mk 
\]%
and hence $K(C)^{k}\leq 1$ and $K(C)\leq 1$ for all uniquely decipherable
codes $C$. This holds for all $k\geq 1.$ Recombining this with (7), letting $%
C$ and $D$ be arbitrary finite uniquely decipherable codes with $C\subseteq 
\overline{D}$, and using now the knowledge that $K(D)\leq 1$, we obtain

\[
K(C)^{k}\leq \tsum\limits_{l=k}^{mk}K(D)^{l}\leq mk\cdot K(D)^{k}\text{ \ \
\ \ \ \ \ \ \ \ \ \ \ \ \ \ \ \ \ \ \ \ \ \ \ \ \ \ \ \ \ (8)} 
\]%
Recall that the definition of $m$ is independent of the choice of $k.$ Thus
(8), being true for all $k\geq 1$, yields the inequality $K(C)\leq K(D):$

\bigskip

\textbf{Extended McMillan Theorem} \textit{If }$C$\textit{\ and }$D$\textit{%
\ are uniquely decipherable codes over the same alphabet, such that every
string in }$C$\textit{\ is a concatenation of strings in }$D$\textit{, then
the Kraft sum of }$C$\textit{\ is less then or equal to the Kraft sum of }$D$%
.

\bigskip

This statement clearly includes the classical McMillan Theorem,
corresponding to the case where $D$ consists of all strings of length $1$.

\bigskip

\bigskip

\bigskip

\textbf{References}

\bigskip

[K] J. Karush, A simple proof of an inequality of McMillan, IRE Trans.
Information\ Theory IT-7 (1961) 118-118

\bigskip

[M] B. McMillan, Two inequalities implied by unique decipherability, IRE
Trans. Information Theory IT-2 (1956) 115-116

\end{document}